\documentclass[12pt]{article}
\usepackage{latexsym}
\usepackage{amssymb}
 
\newcommand{\proof}{{\noindent \bf Proof. }}
\newtheorem{thm}{Theorem}

\newtheorem{prop}{Proposition}

\newcommand{\p}{{\bf p}}
\newcommand{\x}{{\bf x}}
\newcommand{\y}{{\bf y}}

\newcommand{\C}{{\cal C}}

\newcommand{\D}{{\cal D}}
\newcommand{\F}{{\cal F}}
\newcommand{\G}{{\cal G}}

\newcommand{\R}{{\cal P}}

\newcommand{\Y}{{\cal Y}}
\newcommand{\X}{{\cal X}}
\date{}

\begin{document}
\begin{titlepage}
\title{\bf Degree--doubling graph families}

{\author{{\bf  J\'anos K\"orner}
\\e--mail: {\tt korner@di.uniroma1.it} 
\and{\bf Irene Muzi}
\\e--mail:  {\tt irene.muzi@gmail.com}\quad
\medskip
\\Sapienza University of Rome\\
\\ ITALY}}

\maketitle
\begin{abstract}

Let $\G$ be a family of $n$--vertex graphs of uniform degree 2 with the property that the union of any two member graphs has degree four.
We determine the leading term in the asymptotics of the largest cardinality of such a family. Several analogous problems are discussed. 

\end{abstract}
\end{titlepage}

\section{Introduction}

Let $\F$ and $\D$ be two disjoint families of graphs on the same vertex set $[n]=\{1,2,\dots, n\}$. We denote by $M(\F,\D, n)$ the largest cardinality of 
a subfamily $\G\subseteq \F$ with the property that the union of any two of its different members belongs to $\D$. Here the union of two graphs on the same vertex set is the graph whose edge set is the union of those of the two graphs. We are especially interested in the cases when $\D$ is {\em monotone} in the sense 
that if a graph is in $\D$ then any graph containing it also belongs to $\D$. (In our context, a graph contains an other one if they have the same vertex set and their edge sets are in this relation). This framework was introduced in \cite{KMS}. It represents an attempt to describe a consistent part of extremal combinatorics where information theoretic methods seem to 
be relevant. When $\F$ and $\D$ are complementary, then it is clear that without loss of generality, one can suppose that at least one $\G$ of maximum cardinality consists of maximal elements of $\F$. This is true because if we replace a member of  $\G$ by any graph containing it, the union condition 
remains satisfied by the monotonicity of $\D$. Further, if we do this for all the non-maximal member graphs, the family so obtained will have the same cardinality as the original family, since no two members of $\G$ can be contained in a same member of $\F$, for else their union would still be in $\F$, which is impossible by the union condition. Therefore, in this case, a $\G$ of maximum cardinality consists of all the maximal elements of $\F$. Thus if  $\F$ and $\D$ are complementary, our problem reduces to counting the 
maximal elements of the graph family $\F$. Many such enumeration problems have been studied. In several recent papers of this kind information--theoretic methods are used, cf. e. g. \cite{Kahn} dating back such enumeration problems to Dedekind \cite{De} or the often rediscovered Kahn--Lov\'asz theorem \cite{CR}. If, however, $\F$ and $\D$ are far from being complementary, the problems in our framework resemble the graph capacity problem of Claude  Shannon \cite{Sh}.

There are many situations in mathematics where we are given a set of which we have to choose the maximum number of elements any two of which are "more distant'' than a given threshold. In most of these problems distance is measured in terms of a metric. A classical example is the code distance 
problem in information theory \cite{MS}, where in a set of binary sequences of some fixed length we are looking for the largest subset of points any two of which are 
at Hamming distance at least $d$. However, in Shannon's graph capacity problem the elements have to be distant in a {\em structural} instead of a metric sense. For various generalizations and applications of Shannon's problem we refer the reader to the survey article \cite{KO} and the more recent paper \cite{sev}. In the first generalizations of Shannon's problem one considers a set of "distant'' sequences from a finite set, or particular sequences from an infinite set, namely, permutations.
Permutations led the authors of \cite{KMS} to extend the problem of graph capacity to the search for "distant'' Hamilton paths, a rather natural representation of permutations. However, the graph representation naturally leads to graph--theoretic concepts of diversity and  hence problems of an altogether different kind. One of the main objectives in this series of papers has been a (so far not very successful) search for common patterns in the optimal constructions. 

In this paper we are studying the growth rate of $M(\F,\D, n)$ in $n$ in cases of graph families defined in terms of degree conditions on the member graphs.
In the last section we introduce a new graph invariant that generalizes the concept of degree doubling in a different direction. 

Note that logarithms and exponentials are to the base 2. 

\section{Main result}

In our simplest problem $\F$ is the family of connected graphs of uniform degree 2 while $\D$ is the family of all graphs of (maximum) degree at least 4 on the 
common vertex set $[n]$. In other words, the graphs in $\F$ are the Hamilton cycles with vertex set $[n]$. Let us write $Q(n)=M(\F, \D, n)$. We will show that, 
essentially, $Q(n)$ grows like the square root of $n!$. More precisely, we have

\begin{thm}\label{main}
$$ \frac{(n-1)!}{2\cdot\lfloor n/2 \rfloor ! (1+\sqrt{2})^n}\leq Q(n)\leq \frac{n!}{\lfloor n/2 \rfloor ! 2^{\lfloor n/2 \rfloor}}$$
\end{thm}

\proof

We start by proving the upper bound. To this end, suppose for a moment that $n$ is even and let $P$ be a perfect matching, i.e., a graph of uniform 
degree 1, with vertex set $[n]$. Let further $\C=\C(P)$ be the family of all Hamilton cycles containing $P$ as a subgraph. We claim that the union of 
any two graphs in $\C$ has maximum degree strictly less than 4. As a matter of fact, if the union of two graphs has degree 4, then for at least one vertex 
$x \in [n]$ the union has degree 4, meaning that the sets of its incident edges in the two graphs must be disjoint. However, since the two cycles contain a common edge incident to $x$,
the one in the perfect matching $P$, we have a contradiction. This means that if $\G$ is a family of Hamilton cycles in which the union of any two members 
has degree 4,  then $\G$ can contain at most one cycle from $\C$. We claim that
\begin{equation}\label{eq:cnt}
|\C|=\frac{(n/2)!2^{n/2}}{n}.
\end{equation}
To verify this claim note that any linear order of the edges in $P$ combined with any given orientation of the edges of $P$ defines an oriented Hamilton path for 
the vertex set $[n]$. The first vertex of this path is the starting point of the first edge in the linear order, whence the path goes to the endpoint of this edge. 
From here the path continues to the first point of the second edge of $P$, and so on. Since there are $(n/2)!$ orders of the edges in $P$ each of which 
having $2^{n/2}$ orientations of these edges, we see that the number of oriented Hamilton paths on $[n]$ containing a fixed perfect matching is 
$(n/2)!2^{n/2}.$ Every oriented Hamilton path gives rise to an oriented Hamilton cycle in the obvious manner, making its last vertex adjacent to the 
first one. We obtain each oriented Hamilton cycle exactly $n/2$ times in this manner. Further, each Hamilton cycle so generated will appear with both of its 
orientations. In conclusion, every Hamilton cycle containing $P$ is obtained $n$ times which gives (\ref{eq:cnt}). 

It is obvious, by symmetry, that every Hamilton cycle is contained in $\C(P)$ for the same number of perfect matchings $P$ of $[n]$.  Thus, by double counting, considering that the total number of Hamilton cycles is $\frac{(n-1)!}{2}$, we 
obtain that
$$|\G| \leq \frac{n!}{(n/2)!2^{n/2}}.$$
Hence, for $n$ even, 
$$Q(n)\leq \frac{n!}{(n/2)!2^{n/2}}\leq (n/2)!2^{n/2},$$
where for the last inequality we use the obvious bound ${n \choose n/2}\leq 2^n$.

The case of an odd $n$ is similar. We fix a perfect near--matching $P^*$ by which we mean a graph of uniform degree 1 with $n-3$ vertices and a path connecting 
the remaining 3 vertices. As before, we consider the set $\C(P^*)$ of all the Hamilton cycles containing this $P^*$. As in the case of $n$ even, at most one of these cycles can be 
in any family $\G$ satisfying our condition on the degree of pairwise graph unions. The number of Hamilton cycles containing our $P^*$ is now 
$$|\C(P^*)|=\frac{(\lfloor n/2 \rfloor)!2^{\lfloor n/2 \rfloor}}{(n-1)}.$$
Just like in the previous case, we get
$$Q(n) \leq \frac{n!}{(\lfloor n/2 \rfloor)!2^{\lfloor n/2 \rfloor}} $$
as stated in the theorem.
For convenience, we note the following somewhat weaker but nicer form of the bound: 
\begin{equation}\label{eq:fi}
Q(n) \leq  (\lceil n/2\rceil )!2^{\lceil n/2 \rceil}
\end{equation}
for every $n$.

Let us now turn to lower bounding $Q(n)$. To this end, we will use a greedy algorithm to exhibit a large enough family of Hamilton cycles with the required properties. At each step in the algorithm we choose an arbitrary Hamilton cycle and eliminate from the choice space all those incompatible with the chosen one. 
This procedure goes on until the choice space becomes empty. 

To analyze this algorithm, we need an upper bound on the number of the cycles incompatible with a fixed Hamilton cycle $H$. Clearly, a cycle $C$ is incompatible with $H$ if and only if the set of their common edges covers all the vertices in $[n]$. Every such covering contains a minimal covering. In a minimal covering the edges are partitioned into 
single edges and paths of two edges. These are the connected components of the underlying graph. Obviously,  the same covering may contain several minimal coverings. Let us fix a minimal covering and let $s$ be the number of its connected components. Then $\lceil\frac{n}{3}\rceil \leq s\leq \lfloor n/2 \rfloor$ and the number of  adjacent edge pairs in the covering is 
$n-2s$. In consequence, as in the first part of the proof, we see that 
the number of Hamilton cycles whose intersection with $H$ contains our fixed minimal covering is 

$$2^s(s-1)!$$
while the number of minimal coverings with $s$ connected components is 
$${s \choose n-2s}.$$
We conclude that the number of Hamilton cycles that are incompatible with a fixed one is upper bounded by
\begin{equation}\label{eq:su}
\sum_{s=\lceil \frac{n}{3}\rceil}^{\lfloor\frac{n}{2}\rfloor}{s \choose n-2s}2^s(s-1)!.
\end{equation}
Notice that 
$$\sum_{s=\lceil \frac{n}{3}\rceil}^{\lfloor n/2 \rfloor}{s \choose n-2s}2^s(s-1)!\leq 
\left(\sum_{s=\lceil \frac{n}{3}\rceil}^{\lfloor n/2 \rfloor}{n \choose n-2s}2^s\right)\lfloor n/2 \rfloor!$$
Further, rewriting ${n \choose n-2s}2^s={n \choose 2s}(\sqrt{2})^{2s}$ we see that the right--hand side of our last inequality can be further bounded by
$$\left(\sum_{s=1}^n {n \choose s}(\sqrt{2})^{s}\right)\lfloor n/2 \rfloor!=(1+\sqrt{2})^n \lfloor n/2 \rfloor!.$$
This means that the greedy algorithm will eliminate at most $(1+\sqrt{2})^n \lfloor n/2 \rfloor!$ cycles at each step, yielding a cycle family $\G$ with the 
desired union property and containing at least
$$\frac {(n-1)!}{2(1+\sqrt{2})^n \lfloor n/2 \rfloor!}$$
cycles, as claimed for the lower bound.

\hfill$\Box$

Next we turn to the general case of graphs of uniform degree 2. Let therefore $\F$ be the family of graphs of constant degree 2, while as before, $\D$ is the 
family of graphs of maximum degree 4, on the same vertex set $[n]$. We denote
$$R(n)=M(\F, \D, n).$$
Obviously, 
$$Q(n) \leq R(n).$$

However, as we shall see, $R(n)$ grows substantially faster than $Q(n)$. As our next statement shows, the family of all those graphs whose connected components 
are triangles is essentially optimal.

\begin{thm}\label{gen}
$$ c\cdot \frac{n!}{\lfloor n/3 \rfloor !\;6^{n/3}} \leq R(n)\leq e^{\sqrt{n}}\frac{n!}{\lfloor n/3 \rfloor !}$$
where the $c$ in the lower bound is an absolute constant. This constant is 1 if $n$ is a multiple of 3.
\end{thm}

\proof

Let us begin by establishing the lower bound. For simplicity, let us suppose for the time being that $n$ is a multiple of 3. Let $F$ and $G$ be two different graphs with vertex set $[n]$ both of which have only triangles as connected components. We claim that their union contains at least one vertex of degree 4.
In fact, suppose that this is not the case. As we have already established, if the union of two graphs of uniform degree 2 does not have vertices of degree 4, then 
this intersection has no isolated points. This then implies that the intersection contains at least two edges of every triangle of both graphs. However, two edges of a triangle define that triangle. In other words, each triangle of $F$ coincides with some triangle of $G$, and this means that the two graphs coincide; a contradiction. As it is easily seen that if $n$ is a multiple of 3, then the number of those graphs on vertex set $[n]$ whose connected components are 
triangles, is 
$$\frac{n!}{ (n/3) !\;6^{n/3}}.$$
This establishes our lower bound if 3 divides $n$. In the opposite case write $n=3q+r$ where $q, r$ are integers with $3<r<6$ and consider those graphs that contain $q$ connected components on $[3q]$ and which coincide on the fixed set $[n]-[3q]$ of at most 5 vertices. Our previous argument applies on $[3q]$.

To prove an almost matching upper bound, let us choose an arbitrary partition $\p$ of $n$, i. e., a sequence of non-necessarily distinct natural integers $n_i, \; i=1,2,\dots, t$ such that 
$\sum_{i=1}^t n_i=n$. The sequence obtained by any permutation of the indices $i$ is considered an other representation of the same partition $\p$. 
Let us further denote by $\F(\p)$ the family of those graphs of uniform degree 2 on the vertex set $[n]$ that have $t$ connected components with 
vertex sets of cardinality $n_i, \; i=1,2,\dots, t$. We claim that 
\begin{equation}\label{up}
M(\F(\p), \D, n)\leq \frac{n!}{\lfloor n/3 \rfloor !\;6^{n/3}}
\end{equation}
To verify this claim, let $k=k(\p)$ be the number of odd integers among the $n_i$ in $\p$. Consider a graph $P(\p)$ with vertex set $[n]$ and having $k$ connected components that are 3-vertex paths, while the rest are single edges. Then $n-3k$ is even and, by its construction, $P(\p)$ has no isolated vertices. This implies that the subfamily $\C(\p)$ of those graphs in $F(\p)$ which contain $P(\p)$ as a subgraph has no two member graphs with a union of degree 4. Further, by symmetry, we see that constructing such subfamilies for all the different copies of $P(\p)$ on the vertex set $[n]$, we obtain a uniform covering of $F(\p)$.
This yields, as in the proof of the upper bound part of Theorem \ref{main}, 
\begin{equation}\label{pes} 
M(\F(\p), \D, n)\leq \frac{|\F(\p)|}{|\C(\p)|}.
\end{equation}
Hence, in order to obtain a proof of (\ref{up}), we will upper bound the cardinality of $\F(\p)$ and lower bound that of $\C(\p)$.
It is easy to see that 
$$|\F(\p)| \leq \frac{n!}{t!2^t \prod_{i=1}^t n_i}\leq \frac{n!}{t!\prod_{i=1}^t n_i}.$$
On the other hand, let $l$ be such that $3k+2l=n$. Then
$$|\C(\p)|\geq \frac{(k+l)!}{n\cdot t! \prod_{i=1}^t n_i}.$$
Substituting the bounds from the last two inequalities into (\ref{pes}) we obtain
$$M(\F(\p), \D, n)\leq \frac{n\cdot n!}{(k+l)!}$$
Observe that 
$$k+l=k+\frac{n-3k}{2}=\frac{n-k}{2}\geq n/3$$
where the last inequality follows from the obvious relation $k\leq n/3$
which allows us to conclude that
\begin{equation}\label{sub}
M(\F(\p), \D, n)\leq\frac{n\cdot n!}{(n/3)!}.
\end{equation}

On the other hand, we obviously have
$$M(\F, \D, n)\leq \sum_{\p} M(\F(\p), \D, n),$$
where $\p$ runs over the partitions of $n$.
We know from the seminal paper of Hardy and Ramanujan \cite{HR} that the number of the partitions of $n$ is less than $\frac{e^{\sqrt{n}}}{n}$. 
Using this estimate in (\ref{sub}) brings us to our upper bound
$$M(\F, \D, n)\leq e^{\sqrt{n}}\frac{n!}{\lfloor n/3 \rfloor !}.$$

\hfill$\Box$

We conclude this section by a slight variant of the problem about Hamilton cycles. In fact, we ask the same problem for Hamilton paths. Hamilton paths (in the oriented case) are a natural representation of permutations and the present problem area grew out of the problem of permutation capacity \cite{KM}. For the next result let $\F^{H}$ be the set of all (non--oriented) Hamilton paths on the vertex set $[n].$ We will show that up to a linear constant the largest cardinality of a 
set of Hamilton paths with the property that the union of any two has degree 4 is the same as the analogous quantity for Hamilton cycles. More precisely, 
we have

\begin{thm}\label{pat}
$$M(\F, \D, n) \geq M(\F^{H}, \D, n)\geq \frac{2}{n-1} \cdot M(\F, \D, n).$$

\end{thm}

\proof

The first inequality follows from our initial observation about maximal elements of a graph family. In fact, Hamilton cycles are maximal elements in 
the family of connected graphs of degree 2 on $[n].$ 

In order to prove the second inequality, let $\C$ be an optimal family of Hamilton cycles, hence $|\C|=M(\F, \D, n)$. Let further, for any pair of distinct 
vertices $\{a,b\}\in {[n] \choose 2}$ the family $\C(a,b)$ consist of those cycles from $\C$ that contain the edge $\{a,b\}.$ Notice that if we drop the edge 
$\{a,b\}$ from each of the cycles in $\C(a,b)$, the resulting Hamilton paths satisfy our condition, implying that
$|\C(a,b)|\leq M(\F^{H}, \D, n).$
We have
$$ n\cdot M(\F, \D, n)=n\cdot |\C|=\sum_{ \{a,b\}\in [n]} |\C(a,b)| \leq {n \choose 2} M(\F^{H}, \D, n),$$
completing the proof.

\hfill$\Box$

\section{Degree doubling and graph distinguishability}

Our previous problem on degree doubling is a special case of the following. Let $G$ be an arbitrary finite simple graph on $n$ vertices. Without loss of  generality we suppose that its vertex set is $[n]$. Let $F$ be a different graph on $[n]$ but isomorphic to $G$. We will say that $F$ and $G$ are Shannon--distinguishable if there exists a vertex $x \in [n]$ such that its neighborhoods in the two graphs are disjoint. Let us denote by $\nu(G)$ the maximum number of pairwise Shannon--distinguishable copies of $G$. If $G$ is a Hamilton cycle on $[n]$ then $\nu(G)=Q(n)$. If $G$ is a digraph, then we can replace neigborhood with out--neighborhood in the definition. (The out--neigborhood of a vertex is the set of those vertices which are the endpoints of edges starting in the vertex). Clearly, 
$\nu(G)$ is 1 if there are no two vertices in $G$ with disjoint neighborhoods. In case of cycle graphs the determination of $\nu(G)$ is very close to the old question of 
permutation capacity of \cite{KM}. 

More importantly, Shannon's classical problem of graph capacity has a natural formulation in these terms.  
Graph capacity corresponds  to the highest rate at which information can be transmitted over a discrete memoryless stationary channel in an error--free manner. 
In Shannon's information theory a 
channel is modeled by a stochastic matrix $W$. The rows of the matrix are indexed by the elements of a finite alphabet $\X$ and the columns by those of a finite alphabet $\Y$. The element at the crossing of the row of index $x$ and the column $y$, traditionally denoted by $W(y|x)$ is the probability that the transmission of the symbol $x$ results in the reception of the symbol $y$. The repeated use of the channel is characterized by a similar transmission matrix, denoted $W^m$. 
If $m$ symbols are transmitted consecutively, then the transmission of the sequence $\x\in \X^m$ of input symbols results in the reception of $\y\in \Y^m$ with probability $W^m(\y|\x)=\prod_{i=1}^mW(y_i|x_i)$ where $x_i$ is the $i$'th coordinate of the sequence $\x$ and $y_i$ is the $i$'th coordinate of $\y$. 
A code $\C \in \X^m$ of length $m$ is a set of input sequences from $\X^m$ no two of which can result in the same output sequence with positive probability.
The rate of $\C$ is the $m$'th root of its cardinality $|\C|$. (More precisely, in the information theoretic literature the rate is the binary logarithm of this quantity). The supremum of all the code rates, $C(W)$ is the {\it zero--error capacity} of the channel. 
We will say that a code is time--symmetric if the codewords are obtainable from one another by a suitable permutation of the coordinates. It is well--known and easy to see that the supremum of the rates of time--symmetric codes achieves capacity, \cite{CK}, Chapter 11. (Note that in the information--theoretic literature such codes are called fixed composition codes).

Shannon \cite{Sh} observed that the determination of zero--error capacity can be formulated in graph theory. For this purpose one can define a graph $G_W$ with vertex set $\X$ in which two vertices are adjacent if the corresponding input symbols cannot result in the same output symbol with positive probability for both. Correspondingly, one can define the graph $G^m_W=G_{W^m}$.
Then the $m$'root of the largest cardinality of a clique in $G_{W^m}$ is the largest rate of a zero--error code for $m$ uses of the channel. The supremum in $m$ of all these rates is $C(G_W)$, the capacity of the graph $G_W$.
We associate with the matrix 
$W$ a directed graph $G=G(W)$ with vertex set $\X\cup \Y$. We draw an edge from $a\in \X\cup \Y$ to $b \in \X \cup \Y$ if either $W(b|a)>0$ or $a \in \Y$. In the same manner we have a digraph $G^m$ for every $m$ and the matrix $W^m$. Given a probability distribution $P$ on $\X$ we denote by 
$G^m(P)$ the possibly empty graph induced by $G^m$ on the set of those vertices in $\X^m$ in which every element $a\in \X$ appears $mP(a)$ times in the coordinates of the sequences. We denote by $\R(\X)$ the set of all probability distributions on the set $\X$. In conclusion we have
 \begin{prop}\label{prop:Sh}
$$C(G)=\sup_{m}\sup_{P\in \R(\X)}\sqrt[m]{\nu(G^m(P))}.$$
 \end{prop}

\proof

It is sufficient to note that to any two sequences $\x$ and $\y$ in the vertex set of $G^m(P)$ there is a permutation of the coordinates of $\x$ that 
transforms $\x$ in $\y$ and leaves invariant the graph $G^m(P)$. 

\hfill$\Box$

\section{Related open problems}

It seems interesting to ask what happens with our original problem if we replace maximum degree with average degree. More precisely, let $\F$ be as before, the 
family of all Hamilton cycles on $[n]$. Let $\D_\alpha$ be the family of all graphs with average degree at least $\alpha$ for some $\alpha >2$. How does 
$M(\F, \D_\alpha, n)$ depend on $\alpha$?

\end{document}